\documentclass[reqno,a4paper]{amsart}
\usepackage{amssymb}%
\usepackage[hmarginratio=1:1]{geometry}%
\usepackage{ifpdf}
\usepackage{listings}
\ifpdf
 \usepackage[hyperindex]{hyperref}
\else
 \expandafter\ifx\csname dvipdfm\endcsname\relax
 \usepackage[hypertex,hyperindex]{hyperref}
 \else
 \usepackage[dvipdfm,hyperindex]{hyperref}
 \fi
\fi
\allowdisplaybreaks[4]
\theoremstyle{plain}
\newtheorem{thm}{Theorem}[section]
\newtheorem{prop}{Proposition}[section]
\newtheorem{cor}{Corollary}[section]

\theoremstyle{remark}
\newtheorem{rem}{Remark}[section]
\theoremstyle{plain} 

\theoremstyle{definition}

\numberwithin{equation}{section}

\begin{document}
\title[]{Some monotonicity rules for quotient of integrals on time scales}
\author{Zhong-Xuan Mao, Xiao-Yue Du, Jing-Feng Tian*}

\address{Zhong-Xuan Mao \\
Hebei Key Laboratory of Physics and Energy Technology\\
Department of Mathematics and Physics\\
North China Electric Power University \\
Yonghua Street 619, 071003, Baoding, P. R. China}
\email{maozhongxuan000\symbol{64}gmail.com}

\address{Xiao-Yue Du\\
Hebei Key Laboratory of Physics and Energy Technology\\
Department of Mathematics and Physics\\
North China Electric Power University \\
Yonghua Street 619, 071003, Baoding, P. R. China}
\email{duxiaoyue\symbol{64}ncepu.cn}

\address{Jingfeng Tian\\
Hebei Key Laboratory of Physics and Energy Technology\\
Department of Mathematics and Physics\\
North China Electric Power University\\
Yonghua Street 619, 071003, Baoding, P. R. China}
\email{tianjf\symbol{64}ncepu.edu.cn}

\begin{abstract}
  As an efficient mathematical tool, monotonicity rules play an extremely crucial role in the real analysis field. In this paper, we explore some monotonicity rules for quotient of Delta, Nabla and Diamond-Alpha integrals with variable upper limits and parameters on time scales, respectively. Moreover, we consider the monotonicity rules for quotient of the product of multiple Delta integrals with parameters on time scales. Power series is also concerned for being a special case of integral with parameters on time scales.
\end{abstract}
\thanks{*Corresponding author: Jing-Feng Tian, e-mail: tianjf\symbol{64}%
ncepu.edu.cn}
\maketitle

\section{Introduction}
Obtaining the monotonicity of complex quantities through simple quantities is what we always focus on \cite{Qi-JMAA-2021,Qi-RACSAM-2021,Zhao-RACSAM-2020,Zhu-RACSAM-2021,Mao-BBMS-2023}.
Scholars devote themselves to exploring monotonicity rules for quotient of two series or integrals and their applications \cite{Yang-JMAA-2015,Yang-JIA-2016,Zhu-RACSAM-2020,Mao-CRM-2023,Mao-PAMS-2023,Mao-BIMS-2023,Mao-BMMSS-2024,Mao-BBMS-2023-2}.

In 1955, the following critical monotonicity rule was proposed in \cite{Biernacki-AUMCS-1955}.
\begin{thm} \label{thm0-1}
Suppose two convergent real power series $\Psi(s)=\sum_{l=0}^\infty \psi_l s^l$ and $\Phi(s)=\sum_{q=0}^\infty \phi_l s^l$ are defined on $(-r,r)(r>0)$,
then $\Psi(s)/\Phi(s)$ is increasing (decreasing) on $(0,r)$ if $\{\phi_l/\phi_l\}$ is non-constant, increasing (decreasing) and $\phi_l>0$.
\end{thm}

An interesting monotonicity rule for quotient of two integrals was presented by Cheeger, Gromov and Taylor \cite{Cheeger-JDG-1982} in 1982.
\begin{thm} \label{thm0-2}
Suppose functions $\psi, \phi$ are positive and integrable satisfying $\psi/\phi$ is decreasing, then
\begin{equation} \label{R0-2}
\frac{\int_0^s \psi(\upsilon) \textrm{d} \upsilon}{\int_0^s \phi(\upsilon) \textrm{d} \upsilon},
\end{equation}
is decreasing with respect to $s$.
\end{thm}

Recently, Yang, Qian and Chu et al. \cite{Yang-JIA-2017} established the following important monotonicity rule.
\begin{thm} \label{thm0-3}
Suppose function $T\colon[0,r]\to(0,\infty)$ is decreasing and differentiable. If convergent series $\Psi(s)=\sum_{l=l_0}^\infty \psi_l s^l$ and $\Phi(s)=\sum_{l=l_0}^\infty \phi_l s^l$ defined on $[0,r]$ satisfy that $\psi_l/\phi_l$ is non-constant, increasing (decreasing), $\psi_{l_0}/\phi_{l_0} \geq (\leq) 1$ and $\phi_l>0$, then the function
\begin{equation*}
\frac{T(s)+\Psi(s)}{T(s)+\Phi(s)}=\frac{T(s)+\sum_{l=l_0}^\infty \psi_l s^l}{T(s)+\sum_{l=l_0}^\infty \phi_l s^l},
\end{equation*}
is increasing (decreasing) on $[0,r]$.
\end{thm}


In 1988, the time scale was established by Hilger \cite{Hilger-PHD-1988}. It unified the representation of discrete and continuous. Enormous papers on this topic have been published \cite{Mao-BMMSS-2023,Mao-JAAC-2023,Mao-HJMS-2023,Mao-JMI-2021}. We refer readers to \cite{Bohner-Spriner-2016,Bohner-Spriner-2001,Bohner-Spriner-2003,Georgiev-Springer-2018,Martynyuk-Springer-2016,Georgiev-Springer-2016}
for more information.

In this paper, we provide some monotonicity rules on times scales. In section 2, some monotonicity rules for quotient of two integrals with variable upper limits on time scales are provided. In section 3, some monotonicity rules for quotient of two integrals with parameters on time scales are furnished. We conclude the article in section 4.

\section{Monotonicity rules for quotient of two integrals with variable upper limits on time scales}

Delta integral on time scales \cite{Bohner-Spriner-2016} is denoted by
\begin{equation*} 
\int_\alpha^\beta \varphi(\upsilon) \Delta \upsilon, \qquad \alpha,\beta \in \mathbb{T}.
\end{equation*}
Taking $\mathbb{T}=\mathbb{R}$ and $\mathbb{T}=\mathbb{N}$, Delta integral changes into classical integral and summation, respectively.
\begin{equation*}
\int_\alpha^\beta \varphi(\upsilon) \Delta \upsilon=\int_\alpha^\beta \varphi(\upsilon) \textrm{d} \upsilon, \qquad \text{if } \alpha,\beta \in \mathbb{T}=\mathbb{R},
\end{equation*}
\begin{equation*}
\int_\alpha^\beta \varphi(\upsilon) \Delta \upsilon=\sum_{\upsilon=\alpha}^{\beta-1} \varphi(\upsilon), \qquad \text{if } \alpha,\beta\in \mathbb{T}=\mathbb{N}.
\end{equation*}

We first investigate the monotonicity rule for quotient of two Delta integrals with variable upper limits on time scales.
\begin{thm} \label{thm1-1}
If functions $\psi$ and $\phi>0$ are Delta-integrable with $\psi/\phi$ increasing (decreasing) on $[\alpha,\beta]_{\mathbb{T}}$, then
\begin{equation*}
\frac{\int_{\alpha}^s \psi(\upsilon) \Delta \upsilon}{\int_{\alpha}^s \phi(\upsilon) \Delta \upsilon},
\end{equation*}
is increasing (decreasing) for all  $s \in [\alpha,\beta]_{\mathbb{T}}$.
\end{thm}

\begin{proof}
For convenience, we denote
\begin{equation*}
\Phi(s):=\frac{\int_{\alpha}^s \psi(\upsilon) \Delta \upsilon}{\int_{\alpha}^s \phi(\upsilon) \Delta \upsilon}, \quad s\in [\alpha,\beta]_{\mathbb{T}}.
\end{equation*}
A direct computation gives
\begin{equation*}
\Phi^{\Delta}(s)\int_\alpha^s \phi(\upsilon) \Delta \upsilon \int_\alpha^{\sigma(s)} \phi(\upsilon) \Delta \upsilon=\psi(s) \int_\alpha^{s} \phi(\upsilon) \Delta \upsilon- \phi(s) \int_\alpha^{s} \psi(\upsilon) \Delta \upsilon.
\end{equation*}
Noting that
\begin{equation*}
\psi(s) \int_\alpha^{s} \phi(\upsilon) \Delta \upsilon- \phi(s) \int_\alpha^{s} \psi(\upsilon) \Delta \upsilon = \int_\alpha^s \phi(s)\phi(\upsilon)\Big( \frac{\psi(s)}{\phi(s)} - \frac{\psi(\upsilon)}{\phi(\upsilon)} \Big) \Delta \upsilon \geq (\leq) 0,
\end{equation*}
and
\begin{equation*}
\int_\alpha^s \phi(\upsilon) \Delta \upsilon \int_\alpha^{\sigma(s)} \phi(\upsilon) \Delta \upsilon \geq 0,
\end{equation*}
we obtain that $\Phi^{\Delta}(s) \geq (\leq) 0$ holds for all  $s \in [\alpha,\beta]_{\mathbb{T}}$. Thus $\Phi(s)$ is increasing (decreasing) for all $s \in [\alpha,\beta]_{\mathbb{T}}$.
\end{proof}

Theorem \ref{thm1-2} displays that the monotonicity for $(T(s)+ \int_{\alpha}^s \psi(\upsilon) \Delta \upsilon)/(T(s)+ \int_{\alpha}^s \phi(\upsilon) \Delta \upsilon)$ can be determined by $\psi/\phi$ when $T$ meets some conditions.
\begin{thm}\label{thm1-2}
If functions $\psi$ and $\phi>0$ are Delta-integrable with $\psi/\phi$ increasing (decreasing) on $[\alpha,\beta]_{\mathbb{T}}$ and $\psi(\alpha)/\phi(\alpha) \geq (\leq) 1$, function $T\colon[\alpha,\beta]_{\mathbb{T}}\to[0,\infty)$ is  Delta-differential and decreasing, then
\begin{equation*}
\frac{T(s)+ \int_{\alpha}^s \psi(\upsilon) \Delta \upsilon}{T(s)+ \int_{\alpha}^s \phi(\upsilon) \Delta \upsilon},
\end{equation*}
is increasing (decreasing) on $[\alpha,\beta]_{\mathbb{T}}$.
\end{thm}

\begin{proof}
By setting
\begin{equation*}
\Phi(s)=\frac{T(s)+ \int_{\alpha}^s \psi(\upsilon) \Delta \upsilon}{T(s)+ \int_{\alpha}^s \phi(\upsilon) \Delta \upsilon},
\end{equation*}
we obtain
\begin{eqnarray} \label{f1-1}
&& \Phi^{\Delta}(s)\Big( T(s)+ \int_{\alpha}^s \phi(\upsilon) \Delta \upsilon \Big)\Big( T^\sigma(s)+ \int_{\alpha}^{\sigma(s)} \phi(\upsilon) \Delta \upsilon \Big) \nonumber \\
&=& \Big( T^\Delta(s) + \psi(s) \Big) \Big( T(s)+ \int_{\alpha}^s \phi(\upsilon) \Delta \upsilon \Big) - \Big( T(s)+ \int_{\alpha}^s \psi(\upsilon) \Delta \upsilon \Big) \Big( T^\Delta(s) + \phi(s) \Big) \nonumber\\
&=& T^\Delta(s) \Big( \int_\alpha^s \phi(\upsilon) \Delta \upsilon -\int_\alpha^s \psi(\upsilon) \Delta \upsilon \Big) + T(s) ( \psi(s)- \phi(s) ) \nonumber\\
&+& \psi(s) \int_\alpha^s \phi(\upsilon) \Delta \upsilon- \phi(s) \int_\alpha^s \psi(\upsilon) \Delta \upsilon,
\end{eqnarray}
where $T^\sigma(s)$ means $T(\sigma(s))$.

Noting that $\psi/\phi$ increasing (decreasing) on $[\alpha,\beta]_{\mathbb{T}}$ and $\psi(\alpha)/\phi(\alpha) \geq (\leq) 1$, then we obtain
\begin{equation*}
\psi(s)\geq (\leq) \phi(s),\qquad s \in [\alpha,\beta]_{\mathbb{T}}.
\end{equation*}
Therefore, formula (\ref{f1-1}) leads to
\begin{eqnarray*}
&& \Phi^{\Delta}(s)\Big( T(s)+ \int_{\alpha}^s \phi(\upsilon) \Delta \upsilon \Big)\Big( T^\sigma(s)+ \int_{\alpha}^{\sigma(s)} \phi(\upsilon) \Delta \upsilon \Big) \nonumber \\
&=& T^\Delta(s) \Big( \int_\alpha^s \phi(\upsilon) \Delta \upsilon -\int_\alpha^s \psi(\upsilon) \Delta \upsilon \Big) + T(s) ( \psi(s)- \phi(s) ) \nonumber\\
&+& \psi(s) \int_\alpha^s \phi(\upsilon) \Delta \upsilon- \phi(s) \int_\alpha^s \psi(\upsilon) \Delta \upsilon \\
&=& T^\Delta(s) \int_\alpha^s \big( \phi(\upsilon)-\psi(\upsilon) \big) \Delta \upsilon + T(s) ( \psi(s)- \phi(s) ) \nonumber\\
&+& \int_\alpha^s \phi(s)\phi(\upsilon) \Big( \frac{\psi(s)}{\phi(s)} - \frac{\psi(\upsilon)}{\phi(\upsilon)} \Big) \Delta \upsilon \geq (\leq) 0.
\end{eqnarray*}
Thereby we complete the proof.
\end{proof}

Taking $\mathbb{T}=\mathbb{R}$ and $\mathbb{T}=\mathbb{N}$ in Theorem \ref{thm1-2}, we obtain Corollaries \ref{cor1-1} and \ref{cor1-2}, respectively.
\begin{cor} \label{cor1-1}
If functions $\psi$ and $\phi>0$ are integrable with $\psi/\phi$ increasing (decreasing) on $[\alpha,\beta]$ and $\psi(\alpha)/\phi(\alpha) \geq (\leq) 1$, function $T\colon[\alpha,\beta]\to[0,\infty)$ is Delta-differential and decreasing, then
\begin{equation*}
\frac{T(s)+ \int_{\alpha}^s \psi(\upsilon) \textrm{d} \upsilon}{T(s)+ \int_{\alpha}^s \phi(\upsilon) \textrm{d} \upsilon},
\end{equation*}
is increasing (decreasing) for all $s\in[\alpha,\beta]$.
\end{cor}

\begin{cor} \label{cor1-2}
If functions $\psi$ and $\phi>0$ satisfy $\psi/\phi$ increasing (decreasing) on $[\alpha,\beta]_\mathbb{N}$ and $\psi(\alpha)/\phi(\alpha) \geq (\leq) 1$, function $T\colon[\alpha,\beta]_\mathbb{N}\to[0,\infty)$ is Delta-differential and decreasing function, where $\alpha,\beta\in\mathbb{N}$, then
\begin{equation*}
\frac{T(s)+ \sum_{\upsilon=\alpha}^{s-1} \psi(\upsilon)}{T(s)+ \sum_{\upsilon=\alpha}^{s-1} \psi(\upsilon)},
\end{equation*}
is increasing (decreasing) on  $[\alpha,\beta] \cap \mathbb{N}$.
\end{cor}

Now we concern the case that function $T$ is different in the numerator and denominator.
\begin{thm}\label{thm1-3}
If functions $\psi$ and $\phi>0$ are Delta-integral with $\psi/\phi$ increasing (decreasing) on $[\alpha,\beta]_{\mathbb{T}}$ and $T_1$, $T_2$ are positive and Delta-differential with
\begin{equation*}
\frac{\psi(s)}{\int_\alpha^s \psi(\upsilon)\Delta \upsilon}
\geq (\leq) \frac{T_1^\Delta(s)}{T_1(s)}
\geq (\leq) \frac{T_2^\Delta(s)}{T_2(s)}
\geq (\leq) \frac{\phi(s)}{\int_\alpha^s \phi(\upsilon)\Delta \upsilon},
\end{equation*}
for all $s\in[\alpha,\beta]_{\mathbb{T}}$, then
\begin{equation*}
\frac{T_1(s)+ \int_{\alpha}^s \psi(\upsilon) \Delta \upsilon}{T_2(s)+ \int_{\alpha}^s \phi(\upsilon) \Delta \upsilon},
\end{equation*}
is increasing (decreasing) for all $s\in[\alpha,\beta]_{\mathbb{T}}$.
\end{thm}

\begin{proof}
Denote
\begin{equation*}
\Phi(s)=\frac{T_1(s)+ \int_{\alpha}^s \psi(\upsilon) \Delta \upsilon}{T_2(s)+ \int_{\alpha}^s \phi(\upsilon) \Delta \upsilon}.
\end{equation*}
Direct calculation leads to
\begin{eqnarray*}
&& \Phi^\Delta(s) \Big( T_2(s)+ \int_{\alpha}^s \phi(\upsilon) \Delta \upsilon \Big) \Big( T_2(s)+ \int_{\alpha}^s \phi(\upsilon) \Delta \upsilon \Big)^\sigma \\
&=& \Big( T_1^\Delta(s)+ \psi(s) \Big) \Big( T_2(s)+ \int_{\alpha}^s \phi(\upsilon) \Delta \upsilon \Big)-
\Big( T_1(s)+ \int_{\alpha}^s \psi(\upsilon) \Delta \upsilon \Big) \Big( T_2^\Delta(s)+ \phi(s) \Big) \\
&=& T_1^\Delta(s)T_2(s)-T_1(s)T_2^\Delta(s) +T_2(s)\psi(s)-T_2^\Delta(s) \int_{\alpha}^s \psi(\upsilon) \Delta \upsilon \\
&+& T_1^\Delta(s) \int_{\alpha}^s \phi(\upsilon) \Delta \upsilon-T_1(s) \phi(s) + \psi(s)\int_{\alpha}^s \phi(\upsilon) \Delta \upsilon-\phi(s)\int_{\alpha}^s \psi(\upsilon) \Delta \upsilon\\
&=& T_1(s) T_2(s) \Big( \frac{T_1^\Delta(s)}{T_1(s)} - \frac{T_2^\Delta(s)}{T_2(s)} \Big) +
T_2(s) \int_{\alpha}^s \psi(\upsilon) \Delta \upsilon \Big( \frac{\psi(s)}{\int_\alpha^s \psi(\upsilon)\Delta \upsilon} - \frac{T_2^\Delta(s)}{T_2(s)} \Big) \\
&+& T_1(s) \int_{\alpha}^s \phi(\upsilon) \Delta \upsilon \Big( \frac{T_1^\Delta(s)}{T_1(s)} - \frac{\phi(s)}{\int_\alpha^s \phi(\upsilon)\Delta \upsilon} \Big) + \int_\alpha^s \phi(s)\phi(\upsilon) \Big( \frac{\psi(s)}{\phi(s)} - \frac{\psi(\upsilon)}{\phi(\upsilon)} \Big) \Delta \upsilon \\
&\geq& (\leq) 0.
\end{eqnarray*}
\end{proof}

Due to the fact that Delta and Nabla integral are analogical, Theorems \ref{thm1-1}--\ref{thm1-3} are also hold for Nabla integral. For convenience, we omit the proof here.
\begin{thm} \label{nabla1-1}
If functions $\psi$ and $\phi>0$ are Nabla-integrable with $\psi/\phi$ increasing (decreasing) on $[\alpha,\beta]_{\mathbb{T}}$, then
\begin{equation*}
\frac{\int_{\alpha}^s \psi(\upsilon) \nabla \upsilon}{\int_{\alpha}^s \phi(\upsilon) \nabla \upsilon},
\end{equation*}
is increasing (decreasing) for all $s \in [\alpha,\beta]_{\mathbb{T}}$.
\end{thm}

\begin{thm}\label{nabla1-2}
If functions $\psi$ and $\phi>0$ are Nabla-integrable with $\psi/\phi$ increasing (decreasing) on $[\alpha,\beta]_{\mathbb{T}}$ and $\psi(\alpha)/\phi(\alpha) \geq (\leq) 1$, function $T\colon[\alpha,\beta]_\mathbb{T}\to[0,\infty)$ is Nabla-differential and decreasing, then
\begin{equation*}
\frac{T(s)+ \int_{\alpha}^s \psi(\upsilon) \nabla \upsilon}{T(s)+ \int_{\alpha}^s \phi(\upsilon) \nabla \upsilon},
\end{equation*}
is increasing (decreasing) on $[\alpha,\beta]_{\mathbb{T}}$.
\end{thm}

\begin{thm}\label{nabla1-3}
If functions $\psi$ and $\phi>0$ are Nabla-integral with $\psi/\phi$ increasing (decreasing) on $[\alpha,\beta]_{\mathbb{T}}$ and $T_1$, $T_2$ are positive and Nabla-differential with
\begin{equation*}
\frac{\psi(s)}{\int_\alpha^s \psi(\upsilon)\nabla \upsilon}
\geq (\leq) \frac{T_1^\nabla(s)}{T_1(s)}
\geq (\leq) \frac{T_2^\nabla(s)}{T_2(s)}
\geq (\leq) \frac{\phi(s)}{\int_\alpha^s \phi(\upsilon)\nabla \upsilon},
\end{equation*}
for all $s\in[\alpha,\beta]_{\mathbb{T}}$, then
\begin{equation*}
\frac{T_1(s)+ \int_{\alpha}^s \psi(\upsilon) \Delta \upsilon}{T_2(s)+ \int_{\alpha}^s \phi(\upsilon) \Delta \upsilon},
\end{equation*}
is increasing (decreasing) for all $s\in[\alpha,\beta]_{\mathbb{T}}$.
\end{thm}

As a ``weight'' between Delta and Nabla, Sheng, Fadag and Henderson \cite{Sheng-NARWA-7-2006}
defined the Diamond-Alpha derivative and integral as follows
\begin{equation*}
\psi^{\Diamond_\alpha}(\upsilon)=\alpha \psi^\Delta(\upsilon) +(1-\alpha) \psi^\nabla(\upsilon),
\end{equation*}
and
\begin{equation*}
\int_\beta^\gamma \psi(\upsilon) \Diamond_\alpha \upsilon=\alpha \int_\beta^\gamma \psi(\upsilon) \Delta \upsilon +(1-\alpha) \int_\beta^\gamma \psi(\upsilon) \nabla \upsilon.
\end{equation*}

We also consider the monotonicity rule for quotient of two Diamond-Alpha integrals on time scales.
\begin{thm} \label{diamond}
If functions $\psi$ and $\phi>0$ are continuous and Diamond-Alpha integrable with $\psi/\phi$ increasing (decreasing), $\phi$ increasing (decreasing) on $[\beta,\gamma]_{\mathbb{T}}$, then
\begin{equation*}
\frac{\int_{\beta}^s \psi(\upsilon) \Diamond_\alpha \upsilon}{\int_{\beta}^s \phi(\upsilon) \Diamond_\alpha \upsilon},
\end{equation*}
is increasing (decreasing) for all $s\in[\beta,\gamma]_{\mathbb{T}}$.
\end{thm}

\begin{proof}
Setting
\begin{equation*}
\Phi(s)=\frac{\int_{\beta}^s \psi(\upsilon) \Diamond_\alpha \upsilon}{\int_{\beta}^s \phi(\upsilon) \Diamond_\alpha \upsilon}=
\frac{\alpha \int_{\beta}^s \psi(\upsilon) \Delta \upsilon+(1-\alpha) \int_{\beta}^s \psi(\upsilon) \nabla \upsilon}
{\alpha \int_{\beta}^s \phi(\upsilon) \Delta \upsilon+(1-\alpha) \int_{\beta}^s \phi(\upsilon) \nabla \upsilon},
\end{equation*}
then we have
\begin{eqnarray*}
&& \Phi^\Delta(s) \Big( \int_{\beta}^s \phi(\upsilon) \Diamond_\alpha \upsilon \Big)\Big( \int_{\beta}^s \phi(\upsilon) \Diamond_\alpha \upsilon \Big)^\sigma \\
&=& \Big( \alpha \int_{\beta}^s \psi(\upsilon) \Delta \upsilon+(1-\alpha) \int_{\beta}^s \psi(\upsilon) \nabla \upsilon \Big)^\Delta \Big( \alpha \int_{\beta}^s \phi(\upsilon) \Delta \upsilon+(1-\alpha) \int_{\beta}^s \phi(\upsilon) \nabla \upsilon \Big) \\
&-& \Big( \alpha \int_{\beta}^s \psi(\upsilon) \Delta \upsilon+(1-\alpha) \int_{\beta}^s \psi(\upsilon) \nabla \upsilon \Big) \Big( \alpha \int_{\beta}^s \phi(\upsilon) \Delta \upsilon+(1-\alpha) \int_{\beta}^s \phi(\upsilon) \nabla \upsilon \Big)^\Delta \\
&=& \alpha^2 \Big( \psi(s)\int_\beta^s \phi(\upsilon) \Delta \upsilon - \phi(s)\int_\beta^s \psi(\upsilon) \Delta \upsilon \Big) \\
&+& \alpha(1-\alpha) \Big( \psi(s)  \int_\beta^s \phi(\upsilon) \nabla \upsilon- \phi(s)\int_\beta^s \psi(\upsilon) \nabla \upsilon  \Big) \\
&+& \alpha(1-\alpha) \Big( \big( \int_\beta^s \psi(\upsilon) \nabla \upsilon \big)^\Delta \int_\beta^s \phi(\upsilon) \Delta \upsilon- \big( \int_\beta^s \phi(\upsilon) \nabla \upsilon \big)^\Delta \int_\beta^s \psi(\upsilon) \Delta \upsilon  \Big) \\
&+& (1-\alpha)^2 \Big( \big( \int_\beta^s \psi(\upsilon) \nabla \upsilon \big)^\Delta \int_\beta^s \phi(\upsilon) \nabla \upsilon- \big( \int_\beta^s \phi(\upsilon) \nabla \upsilon \big)^\Delta \int_\beta^s \psi(\upsilon) \nabla \upsilon  \Big).
\end{eqnarray*}

Using the the monotonicity of $\psi/\phi$ and Theorem \ref{nabla1-1}, we get
\begin{equation} \label{f2-1}
\psi(s)\int_\beta^s \phi(\upsilon) \Delta \upsilon - \phi(s)\int_\beta^s \psi(\upsilon) \Delta \upsilon = \int_\beta^s \phi(s) \phi(\upsilon) \Big( \frac{\psi(s)}{\phi(s)}- \frac{\psi(\upsilon)}{\phi(\upsilon)} \Big) \Delta \upsilon \geq (\leq) 0,
\end{equation}
\begin{equation} \label{f2-2}
\psi(s)  \int_\beta^s \phi(\upsilon) \nabla \upsilon- \phi(s)\int_\beta^s \psi(\upsilon) \nabla \upsilon = \int_\alpha^s \phi(s) \phi(\upsilon) \Big( \frac{\psi(s)}{\phi(s)}- \frac{\psi(\upsilon)}{\phi(\upsilon)} \Big) \nabla \upsilon \geq (\leq) 0,
\end{equation}
and $\int_{\beta}^s \psi(\upsilon) \nabla \upsilon/\int_{\beta}^s \phi(\upsilon) \nabla \upsilon$ is increasing (decreasing). Then we obtain
\begin{eqnarray} \label{f2-3}
&& \Big( \frac{\int_{\beta}^s \psi(\upsilon) \nabla \upsilon}{\int_{\beta}^s \phi(\upsilon) \nabla \upsilon} \Big)^\Delta \int_{\beta}^s \phi(\upsilon) \nabla \upsilon \int_{\beta}^{\sigma(s)} \phi(\upsilon) \nabla \upsilon \nonumber \\
&=& \Big( \int_{\beta}^s \psi(\upsilon) \nabla \upsilon \Big)^\Delta \int_{\beta}^s \phi(\upsilon) \nabla \upsilon- \int_{\beta}^s \phi(\upsilon) \nabla \upsilon \Big( \int_{\beta}^s \phi(\upsilon) \nabla \upsilon \Big)^\Delta \geq (\leq) 0.
\end{eqnarray}

Based on the hypotheses that $\psi, \phi$ are continuous and $\phi$ increasing (decreasing), we receive
\begin{eqnarray} \label{f2-4}
&& \Big( \int_\beta^s \psi(\upsilon) \nabla \upsilon \Big)^\Delta \int_\beta^s \phi(\upsilon) \Delta \upsilon- \Big( \int_\beta^s \phi(\upsilon) \nabla \upsilon \Big)^\Delta \int_\beta^s \psi(\upsilon) \Delta \upsilon \nonumber \\
&=& \psi(\sigma(s)) \int_\beta^s \phi(\upsilon) \Delta \upsilon- \phi(\rho(s)) \int_\beta^s \psi(\upsilon) \Delta \upsilon \nonumber \\
&=& \int_\beta^s \phi(\upsilon) \phi(\rho(s)) \Big( \frac{\psi(\sigma(s))}{\phi(\rho(s))} - \frac{\psi(\upsilon)}{\phi(\upsilon)} \Big) \Delta \upsilon \nonumber \\
&\geq&(\leq) \int_\beta^s \phi(\upsilon) \phi(\rho(s)) \Big( \frac{\psi(\sigma(s))}{\phi(\sigma(s))} - \frac{\psi(\upsilon)}{\phi(\upsilon)} \Big) \Delta \upsilon \geq(\leq) 0.
\end{eqnarray}

From formulas (\ref{f2-1})--(\ref{f2-4}) and
\begin{equation*}
\int_{\beta}^s \phi(\upsilon) \Diamond_\alpha \upsilon \Big( \int_{\beta}^s \phi(\upsilon) \Diamond_\alpha \upsilon \Big)^\sigma \geq 0,
\end{equation*}
we have $\Phi^\Delta(s) \geq(\leq) 0$. Thereby we complete the proof.
\end{proof}

Especially, if $\alpha=0, 1$, then Theorem \ref{diamond} leads to Theorems \ref{thm1-1} and \ref{nabla1-1}. If $\alpha=\frac{1}{2}, \frac{1}{3}$, then the following corollary holds.
\begin{cor}
If functions $\psi$ and $\phi>0$ are continuous and Diamond-Alpha integral with $\psi/\phi$ increasing (decreasing), $\phi$ increasing (decreasing) on $[\beta,\gamma]_{\mathbb{T}}$, then
\begin{equation*}
\frac{\int_{\beta}^s \psi(\upsilon) \Delta \upsilon+\int_{\beta}^s \psi(\upsilon) \nabla \upsilon}{\int_{\beta}^s \phi(\upsilon) \Delta \upsilon+\int_{\beta}^s \phi(\upsilon) \nabla \upsilon}, \quad \frac{\int_{\beta}^s \psi(\upsilon) \Delta \upsilon+2\int_{\beta}^s \psi(\upsilon) \nabla \upsilon}{\int_{\beta}^s \phi(\upsilon) \Delta \upsilon+2\int_{\beta}^s \phi(\upsilon) \nabla \upsilon},
\end{equation*}
are increasing (decreasing) for all $s\in[\beta,\gamma]_{\mathbb{T}}$.
\end{cor}

\section{Monotonicity rules for quotient of integrals with parameters on time scales}

We investigate the monotonicity rules for
\begin{equation} \label{*}
\frac{\int_\alpha^\beta \Psi(s,\upsilon) \Delta \upsilon}{\int_\alpha^\beta \Phi(s,\upsilon) \Delta \upsilon},
\end{equation}
\begin{equation} \label{**}
\frac{T(s)+\int_\alpha^\beta \Psi(s,\upsilon) \Delta \upsilon}{T(s)+\int_\alpha^\beta \Phi(s,\upsilon) \Delta \upsilon},
\end{equation}
and
\begin{equation*} \label{***}
\frac{\prod_{l=1}^m \int_{\alpha_l}^{\beta_l} \Psi_l(s,\upsilon) \Delta_l \upsilon}{\prod_{l=1}^m \int_{\alpha_i}^{\beta_l} \Phi_l(s,\upsilon) \Delta_l \upsilon},
\end{equation*}
in this section.

Acquiescently, every integral above is convergent and $\Psi, \Phi, \Psi_i, \Phi_i\colon \mathbb{R}\times \mathbb{T} \to \mathbb{R}^+$.

\begin{prop}
The function (\ref{*}) is increasing (decreasing) with respect to $s$ if and only if
\begin{equation} \label{C1}
\int_\alpha^\beta \frac{\partial \Psi(s,\upsilon)}{\partial s} \Delta \upsilon \int_\alpha^\beta \Phi(s,\upsilon) \Delta \upsilon
- \int_\alpha^\beta \Psi(s,\upsilon) \Delta \upsilon \int_\alpha^\beta \frac{\partial \Phi(s,\upsilon)}{\partial s} \Delta \upsilon \geq (\leq) 0.
\end{equation}
\end{prop}

\begin{proof}
Taking derivative to function (\ref{*}) with respect to $s$ leads to
\begin{equation*}
\frac{\partial}{\partial s} \Big( \frac{\int_\alpha^\beta \Psi(s,\upsilon) \Delta \upsilon}{\int_\alpha^\beta \Phi(s,\upsilon) \Delta \upsilon} \Big)
=\frac{\int_\alpha^\beta \frac{\partial \Psi(s,\upsilon)}{\partial s} \Delta \upsilon \int_\alpha^\beta \Phi(s,\upsilon) \Delta \upsilon
- \int_\alpha^\beta \Psi(s,\upsilon) \Delta \upsilon \int_\alpha^\beta \frac{\partial \Phi(s,\upsilon)}{\partial s} \Delta \upsilon
}{\big(\int_\alpha^\beta \Phi(s,\upsilon) \Delta \upsilon\big)^2}.
\end{equation*}
Thus, function (\ref{*}) is increasing (decreasing) if and only if inequality (\ref{C1}) holds.
\end{proof}

\begin{prop}
Suppose function $T\colon[\alpha,\beta]\to(0,\infty)$ is decreasing and differentiable. If
\begin{equation*}
\Psi(s,\upsilon) \geq (\leq) \Phi(s,\upsilon),
\end{equation*}
\begin{equation*}
\frac{\partial \Psi(s,\upsilon)}{\partial s} \geq (\leq) \frac{\partial \Phi(s,\upsilon)}{\partial s},
\end{equation*}
and inequality (\ref{C1})
holds, then function (\ref{**}) is increasing (decreasing) with respect to $s$.
\end{prop}

\begin{proof}
Taking derivative to function (\ref{**}) with respect to $s$ leads to
\begin{eqnarray} \label{CON2}
&& \frac{\partial}{\partial s} \Big( \frac{T(s)+\int_\alpha^\beta \Psi(s,\upsilon) \Delta \upsilon}{T(s)+\int_\alpha^\beta \Phi(s,\upsilon) \Delta \upsilon} \Big) \nonumber \\
&=& \frac{\int_\alpha^\beta \big( T^\prime(s)+\frac{\partial \Psi(s,\upsilon)}{\partial s} \big) \Delta \upsilon \int_\alpha^\beta \big( T(s)+\Phi(s,\upsilon) \big) \Delta \upsilon
- \int_\alpha^\beta \big( T(s)+\Psi(s,\upsilon) \big) \Delta \upsilon \int_\alpha^\beta \big( T^\prime(s) + \frac{\partial \Phi(s,\upsilon)}{\partial s} \big) \Delta \upsilon
}{\big(T(s)+\int_\alpha^\beta \Phi(s,\upsilon) \Delta \upsilon\big)^2} \nonumber \\
&=& \frac{T^\prime(s) \Big(\int_\alpha^\beta \big( \Phi(s,\upsilon)-\Psi(s,\upsilon) \big)  \Delta \upsilon \Big) +T(s) \Big( \int_\alpha^\beta \big( \frac{\partial \Psi(s,\upsilon)}{\partial s} - \frac{\partial \Phi(s,\upsilon)}{\partial s} \big) \Delta \upsilon \Big) }{\big( T(s)+\int_\alpha^\beta \Phi(s,\upsilon) \Delta \upsilon \big)^2} \nonumber \\
&+& \frac{\int_\alpha^\beta \frac{\partial \Psi(s,\upsilon)}{\partial s} \Delta \upsilon \int_\alpha^\beta \Phi(s,\upsilon) \Delta \upsilon
- \int_\alpha^\beta \Psi(s,\upsilon) \Delta \upsilon \int_\alpha^\beta \frac{\partial \Phi(s,\upsilon)}{\partial s} \Delta \upsilon}{\big( T(s)+\int_\alpha^\beta \Phi(s,\upsilon) \Delta \upsilon \big)^2}.
\end{eqnarray}
Since $T$ is non-negative and decreasing, we have $T(s)\geq 0$ and $T^\Delta(s)\leq0$. Together with
\begin{equation*}
\Psi(s,\upsilon) \geq (\leq) \Phi(s,\upsilon),
\end{equation*}
\begin{equation*}
\frac{\partial \Psi(s,\upsilon)}{\partial s} \geq (\leq) \frac{\partial \Phi(s,\upsilon)}{\partial s},
\end{equation*}
and inequality (\ref{C1}), we deduce that $\frac{\partial}{\partial s} \Big( \frac{T(s)+\int_\alpha^\beta \Psi(s,\upsilon) \Delta \upsilon}{T(s)+\int_\alpha^\beta \Phi(s,\upsilon) \Delta \upsilon} \Big) \geq (\leq) 0$.
\end{proof}

\begin{rem}
In \cite{Mao-BMMSS-2024}, Mao and Tian consider the monotonicity of the function
\begin{equation*}
x\mapsto \frac{\int_{t_0}^x \varphi(x,t) w(x,t) \Delta t}{\int_{t_0}^x \psi(x,t) w(x,t) \Delta t}.
\end{equation*}
\end{rem}

Next we consider some concrete cases.

\subsection{Case 1}

Taking $\Psi(s,\upsilon)=\psi(\upsilon) s^\upsilon$ and $\Phi(s,\upsilon)=\phi(\upsilon) s^\upsilon$, where $\psi,\phi$ are positive functions, then function (\ref{*}) transforms into
\begin{equation*}
\frac{\int_\alpha^\beta \psi(\upsilon) s^\upsilon \Delta \upsilon}{\int_\alpha^\beta \phi(\upsilon) s^\upsilon \Delta \upsilon}.
\end{equation*}

\begin{thm} \label{thm2-1}
Suppose $\psi,\phi>0$ are functions defined on $[\alpha,\beta]_\mathbb{T}$. If $\psi/\phi$ is increasing (decreasing), then
\begin{equation} \label{f3-1}
\frac{\int_\alpha^\beta \psi(\upsilon) s^\upsilon \Delta \upsilon}{\int_\alpha^\beta \phi(\upsilon) s^\upsilon \Delta \upsilon},
\end{equation}
is increasing (decreasing).
\end{thm}

\begin{proof}
Taking $\Psi(s,\upsilon)=\psi(\upsilon) s^\upsilon$ and $\Phi(s,\upsilon)=\phi(\upsilon) s^\upsilon$ in the left hand side of inequality (\ref{C1}) yields
\begin{eqnarray}
&& \int_\alpha^\beta \frac{\partial \Psi(s,\upsilon)}{\partial s} \Delta \upsilon \int_\alpha^\beta \Phi(s,\upsilon) \Delta \upsilon
- \int_\alpha^\beta \Psi(s,\upsilon) \Delta \upsilon \int_\alpha^\beta \frac{\partial \Phi(s,\upsilon)}{\partial s} \Delta \upsilon \nonumber\\
&=& \int_\alpha^\beta \upsilon\psi(\upsilon)s^{\upsilon-1} \Delta \upsilon \int_\alpha^\beta \phi(\upsilon) s^\upsilon \Delta \upsilon
- \int_\alpha^\beta \psi(\upsilon) s^\upsilon \Delta \upsilon \int_\alpha^\beta \upsilon\phi(\upsilon)s^{\upsilon-1} \Delta \upsilon \nonumber\\
&=& \frac{1}{2} \Big( \int_\alpha^\beta \int_\alpha^\beta \upsilon\psi(\upsilon) \phi(s) s^{\upsilon+\nu-1} \Delta \nu \Delta \upsilon
- \int_\alpha^\beta \int_\alpha^\beta \upsilon \psi(\nu)\phi(\upsilon)s^{\upsilon+\nu-1} \Delta \nu \Delta \upsilon \nonumber\\
&+& \int_\alpha^\beta \int_\alpha^\beta \nu\psi(\nu) \phi(\upsilon) s^{\upsilon+\nu-1} \Delta \nu \Delta \upsilon
- \int_\alpha^\beta \int_\alpha^\beta \nu \psi(\upsilon)\phi(\nu)s^{\upsilon+\nu-1} \Delta \nu \Delta \upsilon \Big) \nonumber\\
&=& \frac{1}{2} \int_\alpha^\beta \int_\alpha^\beta (\upsilon-\nu) \phi(\nu)\phi(\upsilon) \Big( \frac{\psi(\upsilon)}{\phi(\upsilon)}- \frac{\psi(\nu)}{\phi(\nu)} \Big) s^{\upsilon+\nu-1} \Delta \nu \Delta \upsilon \nonumber\\
&\geq& (\leq) 0. \nonumber
\end{eqnarray}
So the conclusion holds.
\end{proof}

Specially, letting $\mathbb{T}=\mathbb{N}$, $\alpha=0$ and $\beta=\infty$, then function (\ref{f3-1}) transforms into ratios of power series
\begin{equation*}
\frac{\int_0^\infty \psi(\upsilon) s^\upsilon \Delta \upsilon}{\int_0^\infty \phi(\upsilon) s^\upsilon \Delta \upsilon}=\frac{\sum_{\upsilon=0}^\infty \psi(\upsilon) s^\upsilon }{\sum_{\upsilon=0}^\infty \phi(\upsilon) s^\upsilon},
\end{equation*}
namely, Theorem \ref{thm2-1} leads to Theorem \ref{thm0-1}.
Taking $\mathbb{T}=\mathbb{R}$, $\alpha=0$, $\beta=\infty$ and $\mathbb{T}=q^\mathbb{N}(q>1)$, $\alpha=1$, $\beta=\infty$, then the following corollaries hold. This ratio is also researched in \cite{Yang-JMAA-2015,Yang-RJ-2019}.

\begin{cor} \label{cor2-1-1}
Suppose $\psi,\phi>0$ are functions defined on $[0,\infty)$. If $\psi/\phi$ are increasing (decreasing), then
\begin{equation*}
\frac{\int_0^\infty \psi(\upsilon) s^\upsilon \rm{d} \upsilon }{\int_0^\infty \phi(\upsilon) s^\upsilon \rm{d} \upsilon},
\end{equation*}
is increasing (decreasing).
\end{cor}

\begin{cor} \label{cor2-1-2}
Suppose $\psi,\phi>0$ are functions defined on $q^\mathbb{N}(q>1)$. If $\psi/\phi$ are increasing (decreasing), then
\begin{equation*}
\frac{\sum_{\upsilon=0}^\infty (q-1) q^\upsilon \psi(q^\upsilon) s^{q^\upsilon} }{\sum_{\upsilon=0}^\infty (q-1) q^\upsilon \phi(q^\upsilon) s^{q^\upsilon}},
\end{equation*}
is increasing (decreasing).
\end{cor}

\begin{thm}
Suppose $\psi,\phi>0$ are functions defined on $[\alpha,\beta]_\mathbb{T}$. If $\psi/\phi$ are increasing (decreasing) and $\psi(\alpha)/\phi(\alpha) \geq (\leq) 1$, function $T\colon[\alpha,\beta]\to(0,\infty)$ is differential and decreasing function, then
\begin{equation*}
\frac{T(s)+\int_\alpha^\beta \psi(\upsilon) s^\upsilon \Delta \upsilon}{T(s)+\int_\alpha^\beta \phi(\upsilon) s^\upsilon \Delta \upsilon},
\end{equation*}
is increasing (decreasing) on $[\alpha,\beta]$.
\end{thm}
\begin{proof}
Noting that $\psi(s)\geq (\leq) \phi(s)$, formula (\ref{CON2}) leads to
\begin{eqnarray*}
&& \frac{\partial}{\partial s} \Big( \frac{T(s)+\int_\alpha^\beta \psi(\upsilon) s^\upsilon \Delta \upsilon}{T(s)+\int_\alpha^\beta \phi(\upsilon) s^\upsilon \Delta \upsilon}\Big) \\
&=& \frac{\Big( T^\prime(s)+\int_\alpha^\beta \psi(\upsilon) \upsilon s^{\upsilon-1} \Delta \upsilon \Big) \Big( T(s)+\int_\alpha^\beta \phi(\upsilon) s^\upsilon \Delta \upsilon\Big)}{\big( T(s)+\int_\alpha^\beta \phi(\upsilon) s^\upsilon \Delta \upsilon \big)^2} \\
&-&  \frac{\Big( T(s)+\int_\alpha^\beta \psi(\upsilon) s^\upsilon \Delta \upsilon\Big)\Big( T^\prime(s)+\int_\alpha^\beta \phi(\upsilon) \upsilon s^{\upsilon-1} \Delta \upsilon \Big)}{\big( T(s)+\int_\alpha^\beta \phi(\upsilon) s^\upsilon \Delta \upsilon \big)^2}\\
&=& \frac{T^\prime(s) \Big(\int_\alpha^\beta \big(\phi(\upsilon)-\psi(\upsilon)\big) s^\upsilon \Delta \upsilon \Big) +T(s) \Big( \int_\alpha^\beta \big(\psi(\upsilon)-\phi(\upsilon)\big) \upsilon s^{\upsilon-1} \Delta \upsilon \big) }{\big( T(s)+\int_\alpha^\beta \phi(\upsilon) s^\upsilon \Delta \upsilon \big)^2} \\
&+& \frac{\int_\alpha^\beta \upsilon\psi(\upsilon)s^{\upsilon-1} \Delta \upsilon \int_\alpha^\beta \phi(\upsilon) s^\upsilon \Delta \upsilon
- \int_\alpha^\beta \psi(\upsilon) s^\upsilon \Delta \upsilon \int_\alpha^\beta \upsilon\phi(\upsilon)s^{\upsilon-1} \Delta \upsilon}{\big( T(s)+\int_\alpha^\beta \phi(\upsilon) s^\upsilon \Delta \upsilon \big)^2} \\
&\geq& (\leq) 0.
\end{eqnarray*}
Then the conclusion holds.
\end{proof}

Now we consider the monotonic rule for ratios of product of multiple integrals.
\begin{thm} \label{thm2-3}
Suppose functions $\psi_l,\phi_l\colon[\alpha_l,\beta_l]_{\mathbb{T}_l}\to \mathbb{R}$, $l=1,2,\cdots,m$. If $\psi_l/\phi_l$ is increasing (decreasing) and $\phi_l>0$ for all $l=1,2,\cdots,m$, then the function
\begin{equation*}
\frac{\prod_{l=1}^m \int_{\alpha_l}^{\beta_l} \psi_l(\upsilon) s^\upsilon \Delta_l \upsilon}{\prod_{l=1}^m \int_{\alpha_l}^{\beta_l} \phi_l(\upsilon) s^\upsilon \Delta_l \upsilon},
\end{equation*}
is increasing (decreasing).
\end{thm}

\begin{proof}
Noting that
\begin{eqnarray*}
&& \frac{\partial}{\partial s} \Big( \frac{\prod_{l=1}^m \int_{\alpha_l}^{\beta_l} \psi_l(\upsilon) s^\upsilon \Delta_l \upsilon}{\prod_{l=1}^m \int_{\alpha_l}^{\beta_l} \phi_l(\upsilon) s^\upsilon \Delta_l \upsilon} \Big) \big(\prod_{l=1}^m \int_{\alpha_l}^{\beta_l} \phi_l(\upsilon) s^\upsilon \Delta_l \upsilon\big)^2 \\
&=& \frac{\partial}{\partial s}\big(\prod_{l=1}^m \int_{\alpha_l}^{\beta_l} \psi_l(\upsilon) s^\upsilon \Delta_l \upsilon\big) \prod_{l=1}^m \int_{\alpha_l}^{\beta_l} \phi_l(\upsilon) s^\upsilon \Delta_l \upsilon \\
&-& \prod_{l=1}^m \int_{\alpha_l}^{\beta_l} \psi_l(\upsilon) s^\upsilon \Delta_l \upsilon \frac{\partial}{\partial s}\big(\prod_{l=1}^m \int_{\alpha_l}^{\beta_l} \phi_l(\upsilon) s^\upsilon \Delta_l \upsilon\big)\\
&=&  \sum_{l=1}^m \Big(\int_{\alpha_l}^{\beta_l} \psi_l(\upsilon) \upsilon s^{\upsilon-1} \Delta_l \upsilon \prod_{u\neq l} \int_{\alpha_u}^{\beta_u} \psi_u(\upsilon) s^{\upsilon} \Delta_u \upsilon \prod_{u=1}^m \int_{\alpha_u}^{\beta_u} \phi_u(\upsilon) s^{\upsilon} \Delta_u \upsilon \Big)  \\
&-& \sum_{l=1}^m \Big(\int_{\alpha_l}^{\beta_l} \phi_l(\upsilon) \upsilon s^{\upsilon-1} \Delta_l \upsilon \prod_{u\neq l} \int_{\alpha_u}^{\beta_u} \phi_u(\upsilon) s^{\upsilon} \Delta_u \upsilon \prod_{u=1}^m \int_{\alpha_u}^{\beta_u} \psi_u(\upsilon) s^{\upsilon} \Delta_u \upsilon \Big)  \\
&=& \sum_{l=1}^m \Bigg(  \prod_{u\neq l} \Big( \int_{\alpha_u}^{\beta_u} \psi_u(\upsilon) s^{\upsilon} \Delta_u \upsilon \int_{\alpha_u}^{\beta_u} \phi_u(\upsilon) s^{\upsilon} \Delta_u \upsilon \Big) \times \\
&& \Big( \int_{\alpha_l}^{\beta_l} \psi_l(\upsilon) \upsilon s^{\upsilon-1} \Delta_l \upsilon \int_{\alpha_l}^{\beta_l} \phi_l(\upsilon) s^{\upsilon} \Delta_l \upsilon - \int_{\alpha_l}^{\beta_l} \phi_l(\upsilon) \upsilon s^{\upsilon-1} \Delta_l \upsilon \int_{\alpha_l}^{\beta_l} \psi_l(\upsilon) s^{\upsilon} \Delta_l \upsilon \Big) \Bigg)  \\
&=& \sum_{l=1}^m \Bigg(  \prod_{u\neq l} \Big( \int_{\alpha_u}^{\beta_u} \psi_u(\upsilon) s^{\upsilon} \Delta_u \upsilon \int_{\alpha_u}^{\beta_u} \phi_u(\upsilon) s^{\upsilon} \Delta_u \upsilon \Big) \times \\
&& \frac{1}{2} \int_{\alpha_l}^{\beta_l} \int_{\alpha_l}^{\beta_l} (\upsilon-\nu) \phi_l(\upsilon)\phi_l(\nu) \Big( \frac{\psi_l(\upsilon)}{\phi_l(\upsilon)}- \frac{\psi_l(\nu)}{\phi_l(\nu)} \Big) s^{\upsilon+\nu-1} \Delta_l \nu \Delta_l \upsilon \Bigg) \\
&\geq& (\leq) 0,
\end{eqnarray*}
then the conclusion obtains.
\end{proof}

Taking $\mathbb{T}=\mathbb{R}, \mathbb{N}$ and $q^\mathbb{N}$, Theorem \ref{thm2-3} gives the following corollaries, respectively.
\begin{cor}
If functions $\psi_l,\phi_l\colon[0,\infty)\to \mathbb{R}$
satisfy that $\psi_l/\phi_l$ is increasing (decreasing) and $\phi_l>0$ for all $l=1,2,\cdots,m$, then
\begin{equation*}
\frac{\prod_{l=1}^m \int_{0}^{\infty} \psi_l(\upsilon) s^\upsilon \rm{d} \upsilon}{\prod_{l=1}^m \int_{0}^{\infty} \phi_l(\upsilon) s^\upsilon \rm{d} \upsilon},
\end{equation*}
is increasing (decreasing).
\end{cor}

\begin{cor}
If functions
$\psi_l,\phi_l\colon\mathbb{N}\to \mathbb{R}$
satisfy that $\psi_l/\phi_l$ is increasing (decreasing) and $\phi_l>0$ for all $l=1,2,\cdots,m$, then
\begin{equation*}
\frac{\prod_{l=1}^m \Big( \sum_{\upsilon=0}^{\infty} \psi_l(\upsilon) s^\upsilon \Big)}{\prod_{l=1}^m \Big( \sum_{\upsilon=0}^{\infty} \phi_l(\upsilon) s^\upsilon \Big)},
\end{equation*}
is increasing (decreasing).
\end{cor}

\begin{cor}
If functions
$\psi_l,\phi_l\colon q^\mathbb{N}\to \mathbb{R} (q>1, l=1,2,\cdots,m)$
satisfy that $\psi_l/\phi_l$ is increasing (decreasing) and $\phi_l>0$ for all $l=1,2,\cdots,m$, then
\begin{equation*}
\frac{\prod_{l=1}^m \Big( \sum_{\upsilon=0}^\infty (q-1) q^\upsilon \psi_l(q^\upsilon) s^{q^\upsilon} \Big)}{\prod_{l=1}^m \Big(  \sum_{\upsilon=0}^\infty (q-1) q^\upsilon \phi_l(q^\upsilon) s^{q^\upsilon} \Big)},
\end{equation*}
is increasing (decreasing).
\end{cor}

In fact, each time scale can be different, for instance
\begin{equation*}
\frac{\displaystyle\int_{0}^{\infty} \psi_1(\upsilon) s^\upsilon \textrm{d} \upsilon \sum_{\upsilon=0}^{\infty} \psi_2(\upsilon) s^\upsilon \sum_{\upsilon=0}^\infty 2^\upsilon \psi_3(2^\upsilon) s^{2^\upsilon} \sum_{\upsilon=0}^\infty 3^\upsilon \psi_4(3^\upsilon) s^{3^\upsilon}}
{\displaystyle \int_{0}^{\infty} \phi_1(\upsilon) s^\upsilon \textrm{d} \upsilon \sum_{\upsilon=0}^{\infty} \phi_2(\upsilon) s^\upsilon \sum_{\upsilon=0}^\infty 2^\upsilon \phi_3(2^\upsilon) s^{2^\upsilon} \sum_{\upsilon=0}^\infty 3^\upsilon \phi_4(3^\upsilon) s^{3^\upsilon}},
\end{equation*}
is increasing (decreasing) if $\psi_1,\phi_1 : [0,\infty) \to \mathbb{R}^+$, $\psi_2,\phi_2 : \mathbb{N} \to \mathbb{R}^+$, $\psi_3,\phi_3 : 2^\mathbb{N} \to \mathbb{R}^+$, $\psi_4,\phi_4 : 3^\mathbb{N} \to \mathbb{R}^+$ and $\psi_l/\phi_l(l=1,2,3,4)$ is increasing (decreasing).

\subsection{Case 2}
Letting $\Psi(s,\upsilon)=\psi(\upsilon)\varphi(s,\upsilon)$ and $\Phi(s,\upsilon)=\phi(\upsilon)\varphi(s,\upsilon)$, then function (\ref{*}) transforms into
\begin{equation} \label{f3-2-1}
\frac{\int_\alpha^\beta \psi(\upsilon)\varphi(s,\upsilon) \Delta \upsilon}{\int_\alpha^\beta \phi(\upsilon)\varphi(s,\upsilon) \Delta \upsilon}.
\end{equation}

The following theorem gives a monotonicity rule for function \ref{f3-2-1}.
\begin{thm} \label{thm2-4}
Suppose functions $\psi,\phi\colon[\alpha,\beta]_\mathbb{T} \to \mathbb{R}$, $\varphi: \mathbb{R} \times [\alpha,\beta]_\mathbb{T} \to \mathbb{R}^+ $ and $\phi>0$. If
\begin{equation} \label{con2-1}
\frac{\partial}{\Delta \upsilon} \Big( \frac{\psi(\upsilon)}{\phi(\upsilon)} \Big) \geq (\leq) 0,
\qquad
\frac{\partial}{\Delta \upsilon} \Big( \frac{\frac{\partial \varphi(s,\upsilon)}{\partial s}}{\varphi(s,\upsilon)} \Big) \geq (\leq) 0,
\end{equation}
for all $s\in\mathbb{R}$, then
\begin{equation*}
\frac{\int_\alpha^\beta \psi(\upsilon)\varphi(s,\upsilon) \Delta \upsilon}{\int_\alpha^\beta \phi(\upsilon)\varphi(s,\upsilon) \Delta \upsilon},
\end{equation*}
is increasing on $\mathbb{R}$. If
\begin{equation*} \label{con2-2}
\frac{\partial}{\Delta \upsilon} \Big( \frac{\psi(\upsilon)}{\phi(\upsilon)} \Big) \geq (\leq) 0,
\qquad
\frac{\partial}{\Delta \upsilon} \Big( \frac{\frac{\partial \varphi(s,\upsilon)}{\partial s}}{\varphi(s,\upsilon)} \Big) \leq (\geq) 0,
\end{equation*}
for all $s\in\mathbb{R}$, then
\begin{equation*}
\frac{\int_\alpha^\beta \psi(\upsilon)\varphi(s,\upsilon) \Delta \upsilon}{\int_\alpha^\beta \phi(\upsilon)\varphi(s,\upsilon) \Delta \upsilon},
\end{equation*}
is decreasing on $\mathbb{R}$.
\end{thm}

\begin{proof}
We consider the first case, and the second is similar. Based on condition \ref{con2-1}, we obtain
\begin{equation*}
\Big( \frac{\psi(\upsilon)}{\phi(\upsilon)}-\frac{\psi(\nu)}{\phi(\nu)}\Big)\Big( \frac{\frac{\partial \varphi(s,\upsilon)}{\partial s}}{\varphi(s,\upsilon)} - \frac{\frac{\partial \varphi(s,\nu)}{\partial s}}{\varphi(s,\nu)} \Big) \geq 0,
\end{equation*}
for all $\upsilon,\nu\in [\alpha,\beta]_\mathbb{T}$.
Denote
\begin{equation*}
P(s)=\frac{\int_\alpha^\beta \psi(\upsilon)\varphi(s,\upsilon) \Delta \upsilon}{\int_\alpha^\beta \phi(\upsilon)\varphi(s,\upsilon) \Delta \upsilon}, \qquad s\in\mathbb{R}.
\end{equation*}
Direct differentiation gives
\begin{eqnarray*}
&& P^\prime(s) \Big( \int_\alpha^\beta \phi(\upsilon)\varphi(s,\upsilon) \Delta \upsilon \Big)^2\\
&=& \Big( \int_\alpha^\beta \psi(\upsilon)\varphi(s,\upsilon) \Delta \upsilon \Big)^\prime \Big( \int_\alpha^\beta \phi(\upsilon)\varphi(s,\upsilon) \Delta \upsilon \Big) \\
&-& \Big( \int_\alpha^\beta \psi(\upsilon)\varphi(s,\upsilon) \Delta \upsilon \Big) \Big( \int_\alpha^\beta \phi(\upsilon)\varphi(s,\upsilon) \Delta \upsilon \Big)^\prime \\
&=& \Big( \int_\alpha^\beta \psi(\upsilon) \frac{\partial \varphi(s,\upsilon)}{\partial s} \Delta \upsilon \Big) \Big( \int_\alpha^\beta \phi(\upsilon)\varphi(s,\upsilon) \Delta \upsilon \Big) \\
&-& \Big( \int_\alpha^\beta \psi(\upsilon)\varphi(s,\upsilon) \Delta \upsilon \Big) \Big( \int_\alpha^\beta \phi(\upsilon)\frac{\partial \varphi(s,\upsilon)}{\partial s} \Delta \upsilon \Big)\\
&=& \int_\alpha^\beta \int_\alpha^\beta \Big( \psi(\upsilon) \frac{\partial \varphi(s,\upsilon)}{\partial s} \phi(\nu)\varphi(s,\nu) - \psi(\upsilon)\varphi(s,\upsilon)\phi(\nu)\frac{\partial \varphi(s,\nu)}{\partial s} \Big) \Delta \nu \Delta \upsilon\\
&=& \int_\alpha^\beta \int_\alpha^\beta \psi(\upsilon) \phi(\nu) \varphi(s,\nu) \varphi(s,\upsilon) \Big( \frac{\frac{\partial \varphi(s,\upsilon)}{\partial s}}{\varphi(s,\upsilon)} - \frac{\frac{\partial \varphi(s,\nu)}{\partial s}}{\varphi(s,\nu)} \Big) \Delta \nu \Delta \upsilon \\
&=& \frac{1}{2} \int_\alpha^\beta \int_\alpha^\beta \phi(\upsilon) \phi(\nu) \varphi(s,\upsilon) \varphi(s,\nu) \Big( \frac{\psi(\upsilon)}{\phi(\upsilon)}-\frac{\psi(\nu)}{\phi(\nu)}\Big)\Big( \frac{\frac{\partial \varphi(s,\upsilon)}{\partial s}}{\varphi(s,\upsilon)} - \frac{\frac{\partial \varphi(s,\nu)}{\partial s}}{\varphi(s,\nu)} \Big) \Delta \nu \Delta \upsilon \\
&\geq& 0.
\end{eqnarray*}
Thus we complete the proof.
\end{proof}

Some interesting corollaries are established when taking some special time scales in Theorem \ref{thm2-4}.
\begin{cor}
\cite[Lemma 2.7]{Qi-Hal-2020} Suppose functions $\psi,\phi\colon[\alpha,\beta] \to \mathbb{R}$, $\varphi: \mathbb{R} \times [\alpha,\beta] \to \mathbb{R}^+ $ are integrable and $\phi>0$. If both  $\frac{\partial}{\partial \upsilon} \big( \frac{\psi(\upsilon)}{\phi(\upsilon)} \big)$ and
$\frac{\partial}{\partial \upsilon} \big( \frac{\partial \varphi(s,\upsilon)}{\partial s}/\varphi(s,\upsilon) \big)$
are non-negative or non-positive for all $s \in \mathbb{R}, \upsilon\in[\alpha,\beta]$, then
\begin{equation} \label{R}
\frac{\int_\alpha^\beta \psi(\upsilon)\varphi(s,\upsilon) \textrm{d} \upsilon}{\int_\alpha^\beta \phi(\upsilon)\varphi(s,\upsilon) \textrm{d} \upsilon},
\end{equation}
is increasing on $\mathbb{R}$. If one of them is non-negative and the other is non-positive, then function (\ref{R}) is decreasing on $\mathbb{R}$.
\end{cor}

\begin{rem}
This kind monotonicity rule is also appeared in \cite{Mao-CRM-2023,Mao-PAMS-2023}.
\end{rem}

\begin{cor}
Suppose functions $\psi,\phi: \mathbb{N} \to \mathbb{R}$, $\varphi: \mathbb{R} \times \mathbb{N} \to \mathbb{R}^+ $ and $\phi>0$. If both
\begin{equation*}
\frac{\psi(\upsilon+1)}{\phi(\upsilon+1)} - \frac{\psi(\upsilon)}{\phi(\upsilon)},
\end{equation*}
and
\begin{equation*}
\frac{\frac{\partial \varphi(s,\upsilon+1)}{\partial s}}{\varphi(s,\upsilon+1)} - \frac{\frac{\partial \varphi(s,\upsilon)}{\partial s}}{\varphi(s,\upsilon)},
\end{equation*}
are non-negative or non-positive for all $\upsilon\in\mathbb{N}, s\in\mathbb{R}$, then
\begin{equation} \label{R2}
\frac{\sum_{\upsilon=0}^\infty \psi(\upsilon)\varphi(s,\upsilon)}{\sum_{\upsilon=0}^\infty \phi(\upsilon)\varphi(s,\upsilon)},
\end{equation}
is increasing on $\mathbb{R}$. If one of them is non-negative and the other is non-positive, then function (\ref{R2}) is decreasing on $\mathbb{R}$.
\end{cor}

\begin{cor}
Suppose functions $\psi,\phi: q^\mathbb{N} \to \mathbb{R} (q>1)$, $\varphi: \mathbb{R} \times q^\mathbb{N} \to \mathbb{R}^+ $ and $\phi>0$. If both
\begin{equation*}
\frac{\psi(q\upsilon)}{\phi(q\upsilon)} - \frac{\psi(\upsilon)}{\phi(\upsilon)},
\end{equation*}
and
\begin{equation*}
\frac{\frac{\partial \varphi(s,q\upsilon)}{\partial s}}{\varphi(s,q\upsilon)} - \frac{\frac{\partial \varphi(s,\upsilon)}{\partial s}}{\varphi(s,\upsilon)},
\end{equation*}
are non-negative or non-positive for all $(s,\upsilon) \in \mathbb{R} \times q^\mathbb{N}$, then
\begin{equation} \label{R3}
\frac{\sum_{\upsilon=0}^\infty (q-1) \psi(q^\upsilon)\varphi(s,q^\upsilon)}{\sum_{\upsilon=0}^\infty (q-1) \phi(q^\upsilon)\varphi(s,q^\upsilon)},
\end{equation}
is increasing on $\mathbb{R}$. If one of them is non-negative and the other is non-positive, then function (\ref{R3}) is decreasing on $\mathbb{R}$.
\end{cor}

The monotonicity rule for ratios of product of multiple integrals is established as follow.
\begin{thm} \label{thm2-5}
If functions $\psi_l,\phi_l>0$ defined on $[\alpha_l,\beta_l]_\mathbb{T}$ and $\varphi_l\colon\mathbb{R} \times [\alpha_l,\beta_l]_{\mathbb{T}_l} \to (0,\infty)$ satisfy that $\psi_l/\phi_l(l=1,2,\cdots,m)$ is increasing (decreasing) and
\begin{equation*}
\frac{\partial}{\Delta_l \upsilon} \Big( \frac{\frac{\partial \varphi_l(s,\upsilon)}{\partial s}}{\varphi_l(s,\upsilon)} \Big) \geq (\leq) 0, \qquad s\in \mathbb{R}, \quad \upsilon \in [\alpha_l,\beta_l]_{\mathbb{T}_l}, \quad l=1,2,\cdots,m,
\end{equation*}
then the function
\begin{equation} \label{f2-3}
\frac{\prod_{l=1}^m \int_{\alpha_l}^{\beta_l} \psi_l(\upsilon) \varphi_l(s,\upsilon) \Delta_l \upsilon}{\prod_{l=1}^m \int_{\alpha_l}^{\beta_l} \phi_l(\upsilon) \varphi_l(s,\upsilon) \Delta_l \upsilon},
\end{equation}
is increasing on $\mathbb{R}^+$.
\end{thm}

\begin{proof}
It's easy to check that
\begin{eqnarray*}
&& \frac{\partial}{\partial s} \Big( \frac{\prod_{l=1}^m \int_{\alpha_l}^{\beta_l} \psi_l(\upsilon) \varphi_l(s,\upsilon) \Delta_l \upsilon}{\prod_{l=1}^m \int_{\alpha_l}^{\beta_l} \phi_l(\upsilon) \varphi_l(s,\upsilon) \Delta_l \upsilon} \Big) \big(\prod_{l=1}^m \int_{\alpha_l}^{\beta_l} \phi_l(\upsilon) \varphi_l(s,\upsilon) \Delta_l \upsilon \big)^2 \\
&=& \frac{\partial}{\partial s}\big(\prod_{l=1}^m \int_{\alpha_l}^{\beta_l} \psi_l(\upsilon) \varphi_l(s,\upsilon) \Delta_l \upsilon\big) \prod_{l=1}^m \int_{\alpha_l}^{\beta_l} \phi_l(\upsilon) \varphi_l(s,\upsilon) \Delta_l \upsilon \\
&-& \prod_{l=1}^m \int_{\alpha_l}^{\beta_l} \psi_l(\upsilon) \varphi_l(s,\upsilon) \Delta_l \upsilon \frac{\partial}{\partial s}\big(\prod_{l=1}^m \int_{\alpha_l}^{\beta_l} \phi_l(\upsilon) \varphi_l(s,\upsilon) \Delta_l \upsilon\big)\\
&=&  \sum_{l=1}^m \Big(\int_{\alpha_l}^{\beta_l} \psi_l(\upsilon) \frac{\partial \varphi_l(s,\upsilon)}{\partial s} \Delta_l \upsilon \prod_{u\neq l} \int_{\alpha_u}^{\beta_u} \psi_u(\upsilon) \varphi_u(s,\upsilon) \Delta_u \upsilon \prod_{u=1}^m \int_{\alpha_u}^{\beta_u} \phi_u(\upsilon) \varphi_u(s,\upsilon) \Delta_u \upsilon \Big)  \\
&-& \sum_{l=1}^m \Big(\int_{\alpha_l}^{\beta_l} \phi_l(\upsilon) \frac{\partial \varphi_i(s,\upsilon)}{\partial s} \Delta_l \upsilon \prod_{u\neq l} \int_{\alpha_u}^{\beta_u} \phi_u(\upsilon) \varphi_u(s,\upsilon) \Delta_u \upsilon \prod_{u=1}^m \int_{\alpha_u}^{\beta_u} \psi_u(\upsilon) \varphi_u(s,\upsilon) \Delta_u \upsilon \Big)  \\
&=& \sum_{l=1}^m \Bigg(  \prod_{u\neq l} \Big( \int_{\alpha_u}^{\beta_u} \psi_u(\upsilon) \varphi_u(s,\upsilon) \Delta_u \upsilon \int_{\alpha_u}^{\beta_u} \phi_u(\upsilon) \varphi_u(s,\upsilon) \Delta_u \upsilon \Big) \times  \\
&& \Big( \int_{\alpha_l}^{\beta_l} \psi_l(\upsilon) \upsilon \frac{\partial \varphi_l(s,\upsilon)}{\partial s} \Delta_l \upsilon \int_{\alpha_l}^{\beta_l} \phi_l(\upsilon) \varphi_l(s,\upsilon) \Delta_l \upsilon \\
&& - \int_{\alpha_l}^{\beta_l} \phi_l(\upsilon) \frac{\partial \varphi_l(s,\upsilon)}{\partial s} \Delta_l \upsilon \int_{\alpha_l}^{\beta_l} \psi_l(\upsilon) \varphi_l(s,\upsilon) \Delta_l \upsilon \Big) \Bigg)  \\
&=& \sum_{l=1}^m \Bigg(  \prod_{u\neq l} \Big( \int_{\alpha_u}^{\beta_u} \psi_u(\upsilon) \varphi_u(s,\upsilon) \Delta_u \upsilon \int_{\alpha_u}^{\beta_u} \phi_u(\upsilon) \varphi_u(s,\upsilon) \Delta_u \upsilon \Big) \times \\
&& \frac{1}{2} \int_{\alpha_l}^{\beta_l} \int_{\alpha_l}^{\beta_l} \phi_l(\upsilon) \phi_l(\nu) \varphi_l(s,\upsilon) \varphi_l(s,\nu) \Big( \frac{\psi_l(\upsilon)}{\phi_l(\upsilon)}-\frac{\psi_l(\nu)}{\phi_l(\nu)}\Big)\Big( \frac{\frac{\partial \varphi_l(s,\upsilon)}{\partial s}}{\varphi_l(s,\upsilon)} - \frac{\frac{\partial \varphi_l(s,\nu)}{\partial s}}{\varphi_l(s,\nu)} \Big) \Delta_l \nu \Delta_l \upsilon \Bigg) \\
&\geq& (\leq) 0.
\end{eqnarray*}
\end{proof}

Taking $\mathbb{T}=\mathbb{R}, \mathbb{N}$ and $q^\mathbb{N}$, Theorem \ref{thm2-5} leads to the following corollaries, respectively.
\begin{cor}
If functions $\psi_l,\phi_l\colon[0,\infty)\to(0,\infty)$ and $\varphi_l\colon\mathbb{R} \times [0,\infty) \to (0,\infty) (l=1,2,\cdots,m)$ satisfy that $\psi_l/\phi_l(l=1,2,\cdots,m)$ is increasing (decreasing) and
\begin{equation*}
\frac{\partial}{\partial \upsilon} \Big( \frac{\frac{\partial \varphi_l(s,\upsilon)}{\partial s}}{\varphi_l(s,\upsilon)} \Big) \geq (\leq) 0, \qquad s \in \mathbb{R}, \upsilon\in [0,\infty),
\end{equation*}
for all $l=1,2,\cdots,m$, then the function
\begin{equation*}
\frac{\prod_{l=1}^m \int_{0}^{\infty} \psi_l(\upsilon) \varphi_l(s,\upsilon) \rm{d} \upsilon}{\prod_{l=1}^m \int_{0}^{\infty} \phi_l(\upsilon) \varphi_l(s,\upsilon) \rm{d} \upsilon},
\end{equation*}
is increasing on $\mathbb{R}$.
\end{cor}

\begin{cor}
If functions $\psi_l,\phi_l>0$ defined on $\mathbb{N}$ and $\varphi_l\colon\mathbb{R} \times \mathbb{N} \to (0,\infty) $ satisfy that $\psi_l/\phi_l(l=1,2,\cdots,m)$ is increasing (decreasing) and
\begin{equation*}
\frac{\frac{\partial \varphi_l(s,\upsilon+1)}{\partial s}}{\varphi_l(s,\upsilon+1)}-\frac{\frac{\partial \varphi_l(s,\upsilon)}{\partial s}}{\varphi_l(s,\upsilon)} \geq (\leq) 0, \qquad s\in \mathbb{R}, \quad \upsilon \in \mathbb{N}, \quad l=1,2,\cdots,m,
\end{equation*}
then the function
\begin{equation*}
\frac{\prod_{l=1}^m \Big( \sum_{\upsilon=0}^{\infty} \psi_l(\upsilon)  \varphi_l(s,\upsilon) \Big)}{\prod_{l=1}^m \Big( \sum_{\upsilon=0}^{\infty} \phi_l(\upsilon) \varphi_l(s,\upsilon) \Big)},
\end{equation*}
is increasing on $\mathbb{R}$.
\end{cor}

\begin{cor}
If functions $\psi_l,\phi_l>0$ defined on $q^\mathbb{N}(q>1)$ and $\varphi_l\colon\mathbb{R} \times q^\mathbb{N} \to (0,\infty)$ satisfy that $\psi_l/\phi_l(l=1,2,\cdots,m)$ is increasing (decreasing) and
\begin{equation*}
\frac{\frac{\partial \varphi_l(s,qy)}{\partial s}}{\varphi_l(s,qy)}-\frac{\frac{\partial \varphi_l(s,\upsilon)}{\partial s}}{\varphi_l(s,\upsilon)} \geq (\leq) 0, \qquad s\in \mathbb{R}, \quad \upsilon \in q^\mathbb{N}, \quad l=1,2,\cdots,m,
\end{equation*}
then the function
\begin{equation*}
\frac{\prod_{l=1}^m \Big( \sum_{\upsilon=0}^\infty (q-1) q^\upsilon \psi_l(q^\upsilon) \varphi_l(s,q^\upsilon) \Big)}{\prod_{l=1}^m \Big(  \sum_{\upsilon=0}^\infty (q-1) q^\upsilon \phi_l(q^\upsilon) \varphi_l(s,q^\upsilon) \Big)},
\end{equation*}
is increasing on $\mathbb{R}$.
\end{cor}

The following corollary contains different time scales. Denote $\mathbb{T}_1=[0,\infty)$, $\mathbb{T}_2=\mathbb{N}$, $\mathbb{T}_3=2^\mathbb{N}$ and $\mathbb{T}_4=3^\mathbb{N}$.
\begin{cor}
Suppose $\psi_l,\phi_l\colon\mathbb{T}_l \to \mathbb{R}$, $\varphi_l\colon\mathbb{R} \times \mathbb{T}_l \to \mathbb{R}^+$ and $\phi_l>0$ for all $l=1,2,3,4$. If $\psi_l/\phi_l$ is increasing (decreasing) for all $l=1,2,3,4$ and
\begin{equation*}
\frac{\partial}{\partial \upsilon} \Big( \frac{\frac{\partial \varphi_1(s,\upsilon)}{\partial s}}{\varphi_1(s,\upsilon)} \Big) \geq (\leq) 0,
\end{equation*}
for all $s \in \mathbb{R}, \upsilon \in [0,\infty)$, then the function
\begin{equation*}
\frac{\displaystyle\int_{0}^{\infty} \psi_1(\upsilon) \varphi_1(s,\upsilon) \textrm{d} \upsilon \sum_{\upsilon=0}^{\infty} \psi_2(\upsilon) \varphi_2(s,\upsilon) \sum_{\upsilon=0}^\infty 2^\upsilon \psi_3(2^\upsilon) \varphi_3(s,2^\upsilon) \sum_{\upsilon=0}^\infty 3^\upsilon \psi_4(3^\upsilon) \varphi_4(s,3^\upsilon)}
{\displaystyle\int_{0}^{\infty} \phi_1(\upsilon) \varphi_1(s,\upsilon) \textrm{d} \upsilon \sum_{\upsilon=0}^{\infty} \phi_2(\upsilon) \varphi_2(s,\upsilon) \sum_{\upsilon=0}^\infty 2^\upsilon \phi_3(2^\upsilon) \varphi_3(s,2^\upsilon) \sum_{\upsilon=0}^\infty 3^\upsilon \phi_4(3^\upsilon) \varphi_4(s,3^\upsilon)},
\end{equation*}
is increasing on $\mathbb{R}$.
\end{cor}

\begin{rem}
If
\begin{equation*}
\frac{\partial}{\partial \upsilon} \Big( \frac{\frac{\partial \varphi_1(s,\upsilon)}{\partial s}}{\varphi_1(s,\upsilon)} \Big) \geq (\leq) 0,
\end{equation*}
holds for all $s \in \mathbb{R}, \upsilon \in [0,\infty)$, then we also have
\begin{equation*}
\frac{\frac{\partial \varphi_2(s,\upsilon+1)}{\partial s}}{\varphi_2(s,\upsilon+1)}-\frac{\frac{\partial \varphi_2(s,\upsilon)}{\partial s}}{\varphi_2(s,\upsilon)} \geq (\leq) 0, \qquad s \in \mathbb{R}, \upsilon \in \mathbb{N},
\end{equation*}
\begin{equation*}
\frac{\frac{\partial \varphi_3(s,2\upsilon)}{\partial s}}{\varphi_3(s,2\upsilon)}-\frac{\frac{\partial \varphi_3(s,\upsilon)}{\partial s}}{\varphi_3(s,\upsilon)} \geq (\leq) 0, \qquad s \in \mathbb{R}, \upsilon \in 2^\mathbb{N},
\end{equation*}
and
\begin{equation*}
\frac{\frac{\partial \varphi_4(s,3\upsilon)}{\partial s}}{\varphi_4(s,3\upsilon)}-\frac{\frac{\partial \varphi_4(s,\upsilon)}{\partial s}}{\varphi_4(s,\upsilon)} \geq (\leq) 0, \qquad s \in \mathbb{R}, \upsilon \in 3^\mathbb{N}.
\end{equation*}
\end{rem}

\subsection{Case 3}
The generalized ``monomial'' on time scales (see \cite{Georgiev-Springer-2018}) is given by
\begin{equation} \label{h}
h_m(\upsilon,\nu):=
\begin{cases}
1 & \text{ if } m=0, \\
\int_\nu^\upsilon h_{m-1}(\tau,\nu) \Delta \tau & \text{ if } m=1,2,3,\cdots.
\end{cases}
\end{equation}
It's easy to check that \cite[Table3]{Bohner-ITSF-2011}
\begin{equation*}
h_m(\upsilon,\nu)=
\begin{cases}
\displaystyle \frac{(\upsilon-\nu)^m}{m!}, & \text{if } \mathbb{T}=\mathbb{R}, \\
\displaystyle \frac{(\upsilon-\nu)(\upsilon-\nu-1)\cdots(\upsilon-\nu-m+1)}{m!}, & \text{if } \mathbb{T}=\mathbb{Z}, \\
\displaystyle \prod_{u=0}^{m-1}\frac{\upsilon-q^u \nu}{\sum_{l=0}^u q^l}, & \text{if } \mathbb{T}=q^\mathbb{Z}(q>1).
\end{cases}
\end{equation*}

Therefore, $m!h_m(s,0)=s^m$ holds when $\mathbb{T}=\mathbb{R}$.
So we generalize power series as follows
\begin{equation*}
\sum_{u=u_0}^\infty \psi_u u! h_u(s,s_0), \qquad p_u\in \mathbb{R}.
\end{equation*}

In the same way, taking $\Psi(s,\upsilon)=\psi_\upsilon \upsilon! h_\upsilon(s,s_0)$ and $\Phi(s,\upsilon)=\phi_\upsilon \upsilon! h_\upsilon(s,s_0)$, $\alpha=u_0, \beta=\infty$ and $\mathbb{T}=\mathbb{N}$, then function (\ref{*}) changes into
\begin{equation*}
\frac{\sum_{u=u_0}^{\infty} u! \psi_u h_u(s,s_0) }{\sum_{u=u_0}^{\infty} u! \phi_u h_u(s,s_0)},
\end{equation*}

We establish the following theorems which consider the monotonicity for quotient of two generalized series on time scales.
\begin{thm} \label{thm2-6}
Suppose sequence $\psi_u$ and $\phi_u>0$ satisfy that $\Psi(s):=\sum_{u=u_0}^\infty \psi_u u! h_u(s,s_0)$ and $\Phi(s):=\sum_{u=u_0}^\infty \phi_u u! h_u(s,s_0)$ are converging on $[\alpha,\beta]$. If
\begin{equation} \label{Con}
\frac{h_{m+2}(s,s_0)}{h_{m+1}(s,s_0)} \leq \frac{h_{m+1}(s,s_0)}{h_{m}(s,s_0)},
\end{equation}
for all $m=u-1,u,\cdots$ and $\psi_u/\phi_u$ is increasing (decreasing), then the function
\begin{equation} \label{f3-3-1}
\frac{\Psi(s)}{\Phi(s)}=\frac{\sum_{u=u_0}^\infty \psi_u u! h_u(s,s_0)}{\sum_{u=u_0}^\infty \phi_u u! h_u(s,s_0)},
\end{equation}
is increasing (decreasing) on $[\alpha,\beta]$.
\end{thm}
\begin{proof}
Setting
\begin{equation*}
P(s)=\frac{\Psi(s)}{\Phi(s)}=\frac{\sum_{u=u_0}^\infty \psi_u u! h_u(s,s_0)}{\sum_{u=u_0}^\infty \phi_u u! h_u(s,s_0)},
\end{equation*}
we have
\begin{equation*}
J:=P^\Delta \Phi\Phi^\sigma=\Psi^\Delta\Phi-\Psi\Psi^\Delta.
\end{equation*}

Then $J$ leads to
\begin{eqnarray*}
J(s)
&=& \sum_{u=u_0}^\infty \psi_u u! h_{u-1}(s,s_0) \sum_{u=u_0}^\infty \phi_u u! h_{u}(s,s_0) -
\sum_{u=u_0}^\infty \psi_u u! h_{u}(s,s_0) \sum_{u=u_0}^\infty \phi_u u! h_{u-1}(s,s_0) \\
&=& \sum_{i=u_0}^\infty \sum_{j=u_0}^\infty \phi_i \phi_j i!j!h_{i-1}(s,s_0) h_{j}(s,s_0) \big( \frac{\psi_i}{\phi_i}-\frac{\psi_j}{\phi_j} \big)\\
&=& \sum_{i=u_0}^\infty \sum_{j=u_0}^\infty \phi_i \phi_j i!j!h_{j-1}(s,s_0) h_{i}(s,s_0) \big( \frac{\psi_j}{\phi_j}-\frac{\psi_i}{\phi_i} \big) \\
&=& \frac{1}{2} \sum_{i=u_0}^\infty \sum_{j=u_0}^\infty \phi_i \phi_j i!j! \big( \frac{\psi_i}{\phi_i}-\frac{\psi_j}{\phi_j} \big) \Big( h_{i-1}(s,s_0) h_{j}(s,s_0) - h_{j-1}(s,s_0) h_{i}(s,s_0) \Big) \\
&=& \frac{1}{2} \sum_{i=u_0}^\infty \sum_{j=u_0}^\infty \phi_i \phi_j i!j! h_{i-1}(s,s_0)h_{j-1}(s,s_0)\big( \frac{\psi_i}{\phi_i}-\frac{\psi_j}{\phi_j} \big) \Big( \frac{h_{j}(s,s_0)}{h_{j-1}(s,s_0)} - \frac{h_{i}(s,s_0)}{h_{i-1}(s,s_0)} \Big) \\
&\geq& (\leq) 0.
\end{eqnarray*}

From $\Phi(s)\Phi^\sigma(s)>0$ and $J\geq (\leq) 0$, we get that $P$ is increasing (decreasing). Thereby we complete the proof.
\end{proof}

Based on the definition of $h$, condition (\ref{Con}) holds for $\mathbb{T}=\mathbb{R}$, $\mathbb{T}=\mathbb{N}$ and $\mathbb{T}=q^\mathbb{N}(q>1)$.
In fact, if $\mathbb{T}=\mathbb{R}$, then
\begin{eqnarray*}
&& \frac{h_{m+2}(s,s_0)}{h_{m+1}(s,s_0)} - \frac{h_{m+1}(s,s_0)}{h_{m}(s,s_0)} = \frac{(s-s_0)^{m+2}/(m+2)!}{(s-s_0)^{m+1}/(m+1)!} - \frac{(s-s_0)^{m+1}/(m+1)!}{(s-s_0)^{m}/m!} \\
&=& (s-s_0) \Big( \frac{1}{m+2}-\frac{1}{m+1} \Big) \leq 0.
\end{eqnarray*}
If $\mathbb{T}=\mathbb{N}$, then
\begin{eqnarray*}
&& \frac{h_{m+2}(s,s_0)}{h_{m+1}(s,s_0)} - \frac{h_{m+1}(s,s_0)}{h_{m}(s,s_0)} = \frac{(s-s_0)(s-s_0-1)\cdots(s-s_0-m-1)/(m+2)!} {(s-s_0)(s-s_0-1)\cdots(s-s_0-m)/(m+1)!} \\
&-& \frac{(s-s_0)(s-s_0-1)\cdots(s-s_0-m)/(m+1)!} {(s-s_0)(s-s_0-1)\cdots(s-s_0-m+1)/m!} \\
&=& \frac{s-s_0-m-1}{m+2}-\frac{s-s_0-m}{m+1} = \frac{-s+s_0-1}{m^2+3 m+2} \leq 0.
\end{eqnarray*}
If $\mathbb{T}=q^\mathbb{N}(q>1)$, then
\begin{eqnarray*}
&& \frac{h_{m+2}(s,s_0)}{h_{m+1}(s,s_0)} - \frac{h_{m+1}(s,s_0)}{h_{m}(s,s_0)}
= \frac{\prod_{u=0}^{m+1}\frac{s-q^u s_0}{\sum_{l=0}^uq^l}} {\prod_{u=0}^{m}\frac{s-q^ux_0}{\sum_{l=0}^u q^l}} - \frac{\prod_{u=0}^{m}\frac{s-q^ux_0}{\sum_{l=0}^u q^l}} {\prod_{u=0}^{m-1}\frac{s-q^ux_0}{\sum_{l=0}^uq^l}} \\
&=& \frac{s-q^{m+1} s_0}{\sum_{l=1}^{m+1} q^l}-\frac{s-q^{m} s_0}{\sum_{l=1}^{m} q^l} = \frac{q^{m+1} (s_0-s)}{\sum_{l=1}^{m+1} q^l \sum_{l=1}^{m} q^l} \leq 0.
\end{eqnarray*}

Under these three time scales, function (\ref{f3-3-1}) becomes
\begin{equation*}
\frac{\sum_{u=u_0}^\infty \psi_u (s-s_0)^u}{\sum_{u=u_0}^\infty \phi_u (s-s_0)^u},
\end{equation*}
\begin{equation*}
\frac{\sum_{u=u_0}^\infty \psi_u (s-s_0)(s-s_0-1)\cdots(s-s_0-m+1)}{\sum_{u=u_0}^\infty \phi_u (s-s_0)(s-s_0-1)\cdots(s-s_0-m+1)},
\end{equation*}
and
\begin{equation*}
\frac{\sum_{u=u_0}^\infty \psi_u u! \prod_{u=0}^{m-1}\frac{s-q^ux_0}{\sum_{l=0}^uq^l}}{\sum_{u=u_0}^\infty \phi_u u! \prod_{u=0}^{m-1}\frac{s-q^ux_0}{\sum_{l=0}^uq^l}}.
\end{equation*}

Hence, we have the following three corollaries.
\begin{cor}[Generalized form of Theorem \ref{thm0-1}]
Suppose sequence $\psi_u$ and $\phi_u>0$ satisfy $\Psi(s):=\sum_{u=u_0}^\infty \psi_u (s-s_0)^u$ and $\Phi(s):=\sum_{u=u_0}^\infty \phi_u (s-s_0)^u$ are converging on $[\alpha,\beta]$. Then the function
$
\Psi(s)/\Phi(s)
$
is increasing (decreasing) on $[\alpha,\beta]$.
\end{cor}

\begin{cor}
Suppose sequence $\psi_u$ and $\phi_u>0$ satisfy that
\begin{equation*}
\Psi(s):=\sum_{u=u_0}^\infty \psi_u (s-s_0)(s-s_0-1)\cdots(s-s_0-m+1),
\end{equation*}
and
\begin{equation*}
\Phi(s):=\sum_{u=u_0}^\infty \phi_u (s-s_0)(s-s_0-1)\cdots(s-s_0-m+1),
\end{equation*}
are converging on $[\alpha,\beta]$. Then the function $\Psi(s)/\Phi(s)$ is increasing (decreasing) on $[\alpha,\beta]$.
\end{cor}

\begin{cor}
Suppose sequence $\psi_m$ and $\phi_m>0$ satisfy that
\begin{equation*}
\Psi(s):=\sum_{m=m_0}^\infty \psi_m m! \prod_{u=0}^{m-1}\frac{s-q^u s_0}{\sum_{l=0}^u q^l},
\end{equation*}
and
\begin{equation*}
\Phi(s):=\sum_{m=m_0}^\infty \phi_m m! \prod_{u=0}^{m-1}\frac{s-q^u s_0}{\sum_{l=0}^u q^l},
\end{equation*}
are converging on $[\alpha,\beta]$. Then the function $\Psi(s)/\Phi(s)$ is increasing (decreasing) on $[\alpha,\beta]$.
\end{cor}

However, whether condition (\ref{Con}) holds for arbitrary time scale needs further research.


\end{document}